\documentclass[10pt,twocolumn]{article}
\usepackage[margin=0.75in]{geometry}

\usepackage{cite}
\usepackage{amsmath,amssymb,amsfonts}
\usepackage{mathtools}
\usepackage{algorithmic}
\usepackage{graphicx}
\usepackage{textcomp}
\usepackage{algorithm}
\usepackage{tikz}
\usepackage{ulem}
\usepackage[table]{xcolor}

\usepackage{hyperref}
\usetikzlibrary{positioning,shapes,arrows.meta,patterns,calc}

\newcommand{\Z}{\mathbb{Z}}
\newcommand{\R}{\mathbb{R}}

\def\BibTeX{{\rm B\kern-.05em{\sc i\kern-.025em b}\kern-.08em
    T\kern-.1667em\lower.7ex\hbox{E}\kern-.125emX}}
    
\begin{document}

\title{ILPU: Iterative Laplace-Based Phase Unwrapping via Bi-Level Optimization}

\author{Michael Fedders, Jakob Schattenfroh, Yanglei Wu, Noah Jaitner, Tom Meyer, Jakob Jordan, \\ Jing Guo, Ingolf Sack, and Hossein S. Aghamiry 
\thanks{Funding by German Research Foundation (DFG) projects 513752256 FOR5628, CRC1340, GRK2260 BIOQIC, CRC 1540 EBM, GU 172614-1, German Federation of Industrial Research Associations (AIF) project KK5611902 BM4, and Medical Scientist Pilot Program of the Charité is greatly appreciated. These sponsors had no role in the study design, collection, analysis and interpretation of data, writing of this manuscript or the decision to submit the article for publication.}
\thanks{All the authors are with  Department of Radiology, Charité – Universitatsmedizin Berlin, 10117, Berlin, Germany (email:hossein.aghamiry@charite.de)}}

\maketitle
\begin{abstract}
Phase unwrapping is an essential preprocessing step for phase-based MRI applications including susceptibility mapping, field mapping, thermometry, and MR elastography. We present Iterative Laplace-Based Phase Unwrapping (ILPU), a bi-level optimization algorithm in which a lower-level solver recovers a continuous phase increment from an incremental Poisson equation via the discrete cosine transform (DCT), while an upper-level solver refines an integer offset map through quality-guided spatial regularization and a restricted local search. This coupling yields robust unwrapping in low-SNR regions through adaptive smoothness penalties and quality-weighted regularization. We evaluated ILPU on 2D and 3D brain MRI phase images against manually unwrapped reference data, using standard Laplace unwrapping, Flynn, and SEGUE as comparison methods. In 2D, ILPU achieves accuracy comparable to SEGUE. In 3D, ILPU attains a relative error of 2.12\% versus 67.59\% for SEGUE and 81.02\% for Laplace, demonstrating a clear advantage in volumetric unwrapping. The algorithm operates at $O(N\log N)$ complexity per iteration through DCT-based Laplacian estimation, and is numerically faster than both Flynn and SEGUE while preserving superior accuracy. These results indicate that the bi-level optimization framework provides a robust and computationally efficient solution for phase unwrapping in MRI.
\end{abstract}

\noindent\textbf{Keywords:} Phase unwrapping, Bi-level optimization, Laplacian unwrapping, Magnetic resonance imaging, Elastography.
\vspace{0.5em}

\section{Introduction}
\label{sec:introduction}

Phase unwrapping recovers the continuous phase field from wrapped measurements confined to the principal value range, serving as an essential preprocessing step for phase-based magnetic resonance imaging (MRI) applications. The complex MRI signal encodes tissue properties in both magnitude and phase, but the arctangent operation used to extract the phases inherently wraps values to $(-\pi, \pi]$, creating discontinuities where the true phase exceeds these bounds. Accurate phase unwrapping enables numerous clinical applications including quantitative susceptibility mapping (QSM) for tissue characterization and iron quantification in neurodegenerative diseases \cite{langkammer2012qsm}, proton resonance frequency shift-based thermometry for monitoring thermal ablation procedures \cite{rieke2008mr}, B$_0$ field mapping for geometric distortion correction in echo-planar imaging \cite{jezzard1995correction}, phase-contrast flow imaging for cardiovascular assessment \cite{markl20124d}, and MR elastography (MRE) for tissue mechanical property characterization \cite{muthupillai1995mre,sack2015mre,sack2023magnetic}. MRE encodes tissue displacement from externally applied mechanical vibrations in the phase images, where unwrapping errors propagate directly into stiffness estimates, making phase quality critical across a broad range of clinical applications, including liver fibrosis staging and tumor characterization. The challenges are particularly important in MRE due to low displacement amplitudes in stiff tissues, discontinuities at organ boundaries, and the need for robust unwrapping across multiple vibration frequencies and motion encoding directions \cite{ma2023mr}.

Despite decades of research, robust phase unwrapping in MRI remains challenging due to several fundamental difficulties. Phase measurements in low signal-to-noise ratio (SNR) regions become unreliable, with noise amplified through the arctangent function generating residues, i.e., pixels where closed-loop path integration yields non-zero multiples of $2\pi$ \cite{schweser2016foundations}. These residues violate the integrability conditions required for unique unwrapping solutions and can propagate errors throughout the image. Near air-tissue interfaces such as paranasal sinuses and auditory canals, large magnetic susceptibility differences create steep phase gradients that often exceed the Nyquist sampling criterion, resulting in spatial aliasing manifested as multiple $2\pi$ wraps between adjacent pixels \cite{witoszynskyj2009phase}. 

Existing phase unwrapping algorithms can be broadly categorized into local and global methods. Local methods rely on local consistency and integration paths. They typically start from high-quality seed points or phase-coherent regions and propagate the unwrapping decision to neighboring pixels or voxels. Because the decision is made locally and then propagated, these methods are usually fast and effective when a reliable quality map is available, but they are sensitive to the starting seed, path selection, and local error propagation. This family includes classical path-following and quality-guided methods such as those of Goldstein et al.~\cite{goldstein1988satellite} and Herr{\'a}ez et al.~\cite{herraez2002fast}, as well as regional methods such as PRELUDE and SEGUE (Speedy rEgion-Growing algorithm for Unwrapping Estimated phase)~\cite{karsa2019segue}. PRELUDE, a gold standard among region-growing methods, groups voxels into phase-coherent regions and merges them by integer offsets. SEGUE builds on PRELUDE and addresses its speed limitation, and serves as one of the reference methods in this work.

Global methods enforce consistency across the entire image or volume, typically through energy minimization or global constraint satisfaction, and can be broadly categorized as either discrete or continuous. Discrete global methods optimize integer offsets or discontinuities directly. This includes network flow methods by Costantini~\cite{costantini1998novel} and Chen and Zebker~\cite{chen2000network,chen2001two}, as well as statistical formulations based on Bayesian inference and Markov random fields~\cite{nico2000bayesian,ying2006unwrapping}. 

Continuous global methods estimate the phase field by minimizing a continuous objective over the full image. This includes the minimum $L^p$-norm formulation by Ghiglia and Romero~\cite{ghiglia1996minimum}, weighted least-squares methods, and Laplacian or Poisson-based approaches~\cite{ghiglia1994robust,schofield2003fast,volkov2003deterministic,li2014integrated}. Laplacian approaches exploit the fact that, under appropriate regularity conditions, the Laplacian of the unwrapped phase can be inferred directly from the wrapped observations. This enables a fast recovery of the phase through solution of a Poisson equation with suitable boundary conditions, often with discrete cosine transform (DCT) solvers and $O(N \log N)$ complexity~\cite{schofield2003fast}. Weighted variants incorporate quality maps by assigning reliability-based weights to phase gradients, but pure Laplacian approaches still lack explicit mechanisms to resolve integer offset ambiguities that may persist in noisy regions or in the presence of rapid phase variations. The Flynn method~\cite{flynn1997minimum} also belongs to this group. It formulates unwrapping as a global optimization problem that minimizes weighted discontinuities, using quality maps to concentrate unavoidable errors in low-reliability regions. Although this method is computationally less convenient for large volumetric datasets, Flynn remains an accurate and well-established reference for 2D data. For this reason, we include it as the other reference throughout this work.

Bi-level optimization provides a natural framework for hierarchical problems where one optimization task (the upper level or leader) constrains another (the lower level or follower) \cite{colson2007overview}. This structure arises naturally in inverse problems where lower-level reconstruction parameters are optimized by upper-level criteria and vice versa. 
For phase unwrapping, bi-level optimization enables natural decomposition. 

This paper presents ILPU (Iterative Laplace-Based Phase Unwrapping), a bi-level optimization algorithm for MRI phase unwrapping. The lower level employs DCT-based Laplacian unwrapping to estimate regularized phase increments with $O(N \log N)$ complexity, while the upper level refines integer-valued offset maps through quality-guided spatial regularization using a restricted local search. Robust performance in low-SNR regions is achieved through adaptive regularization scheduling and quality-guided residual smoothing.

The phase wrapping challenges are particularly notable in MRE due to low displacement amplitudes in stiff tissues, discontinuities at organ boundaries, and the need for robust unwrapping across multiple vibration frequencies and spatial
directions \cite{sack2023magnetic}. Several MRE-specific unwrapping strategies have been proposed~\cite{pmid37030683,pmid25485426,pmid24942537,pmid21666289}, yet none combines the accuracy and computational efficiency needed for routine post-processing of large datasets. In this work we use six-dimensional MRE (6D-MRE) data~\cite{sack2023magnetic}, comprising three spatial dimensions, time-resolved wave motion, motion-encoding directions, and vibration frequencies. This data structure simultaneously captures several unwrapping challenges present in the field, including inter-slice inconsistencies, direction-dependent artifacts, and amplitude overshoots. We focus primarily on 2D unwrapping because many MRE acquisition protocols remain two-dimensional~\cite{pmid40720631,pmid37705467}, while also evaluating the method on 3D data~\cite{pmid35691940}.

\section{Methodology}
\label{sec:methodology}

\subsection{Problem Formulation}
\label{sec:problem_formulation}

Let $\phi \in \R^{N}$ denote the unwrapped phase field defined on a discrete 3D grid with $N = H \times W \times S$ voxels, where $H$, $W$, and $S \in \mathbb{N}$ are the image height, width, and number of slices, respectively. For 2D data, $S = 1$ and $N = H \times W$. The wrapped measurement $\phi_w \in (-\pi, \pi]^N$ relates to the true unwrapped phase $\phi$ through the phase-jump model
\begin{equation}
\phi = \phi_w + 2\pi \mathbf{k},
\label{eq:unwrap_relation}
\end{equation}
where $\mathbf{k} \in \Z^N$ is the integer-valued offset field. The unwrapping problem seeks to recover $\mathbf{k}$ and therefore $\phi$ from $\phi_w$.

Following Itoh~\cite{itoh1982analysis} and Juarez-Salazar et al.~\cite{juarezsalazar2014phase}, if the true phase satisfies the Nyquist sampling condition, i.e., $\|\mathbf{D}\phi\|_\infty < \pi$, the gradients of the wrapped and unwrapped phases are related by
\begin{equation}
\mathbf{D}\phi = \mathcal{W}[\mathbf{D}\phi_w],
\label{eq:itoh}
\end{equation}
where $\mathbf{D} \in \R^{dN \times N}$ stacks first-order finite differences along the active spatial directions, with $d=2$ for 2D data and $d=3$ for 3D data, and $\mathcal{W}[\cdot]$ is the element-wise wrapping operator to $(-\pi, \pi]$. Equation~\eqref{eq:itoh} holds because the differences of the integer-offset term $2\pi\mathbf{k}$ in~\eqref{eq:unwrap_relation} are integer multiples of $2\pi$, which are removed by $\mathcal{W}$.

Global phase-unwrapping methods recover $\phi$ as the minimizer of
\begin{equation}
\min_{\phi} \; \|\mathbf{D}\phi - \mathcal{W}[\mathbf{D}\phi_w]\|_p,
\label{eq:global_unwrap}
\end{equation}
for some norm $p$~\cite{ghiglia1996minimum}. For $p = 1$, Eq.~\eqref{eq:global_unwrap} belongs to the class of minimum-discontinuity formulations represented by Flynn's method~\cite{flynn1997minimum}. For $p = 2$, minimizing the squared norm gives the least-squares Poisson equation
\begin{equation}
\mathbf{L}\phi = \rho, \qquad
\rho = \mathbf{D}^{\top}\mathcal{W}[\mathbf{D}\phi_w],
\label{eq:poisson}
\end{equation}
where $\mathbf{L} = \mathbf{D}^{\top}\mathbf{D}$ is the discrete Laplacian. Equation~\eqref{eq:poisson} can be solved efficiently in the DCT domain with $O(N\log N)$ complexity~\cite{ghiglia1994robust,schofield2003fast}.

Although the $\ell_2$ formulation~\eqref{eq:poisson} is computationally efficient, it suffers from two well-documented limitations~\cite{juarezsalazar2014phase,ghiglia1996minimum}. First, \textit{residue spreading}: because the minimization is global, phase residues caused by noise or undersampling are distributed across the image rather than confined to their originating voxels. Second, \textit{over-smoothing}: the least-squares solution may attenuate genuine high-frequency phase variations near tissue boundaries and in regions of rapid spatial phase variation.

Rather than solving Eq.~\eqref{eq:poisson} from scratch, we update the phase estimate incrementally \cite{gholamiaghamiry2017}. Let $\phi^{(t)}$ denote the current continuous estimate at iteration $t$. The next estimate is 
\begin{equation}
\phi^{(t+1)} = \phi^{(t)} + \delta\phi^{(t+1)},
\label{eq:increment_definition}
\end{equation}
where $\delta\phi^{(t+1)}$ is the phase increment. Subtracting $\phi^{(t)}$ from Eq.~\eqref{eq:unwrap_relation} gives
\begin{equation}
\delta\phi^{(t+1)} = \phi_w + 2\pi\mathbf{k}^{(t)} - \phi^{(t)}.
\label{eq:increment_model}
\end{equation}
If $\|\mathbf{D}\delta\phi\|_\infty < \pi$, the incremental Itoh relation is
\begin{equation}
\mathbf{D}\delta\phi^{(t+1)}
= \mathcal{W}[\mathbf{D}\mathcal{W}[\phi_w + 2\pi\mathbf{k}^{(t)} - \phi^{(t)}]].
\label{eq:increment_itoh}
\end{equation}

In practice, the residual phase contains noise and local inconsistencies. For a candidate offset map $\mathbf{k}$, we therefore apply a quality-guided smoothing operator $\mathcal{S}$ to the residual as
\begin{equation}
\Delta^{(t)}(\mathbf{k})
= \mathcal{W}[\mathcal{S}[\phi_w + 2\pi\mathbf{k} - \phi^{(t)}]].
\label{eq:regularized_residual_main}
\end{equation}
In high-quality regions, $\mathcal{S}$ approaches the identity operator, whereas in low-quality regions it suppresses noise-driven fluctuations. Using the regularized residual, the phase increment is estimated in the least-squares sense:
\begin{equation}
\delta\phi^{(t+1)}(\mathbf{k})
\in \arg\min_{\delta\phi \in \R^N}
\|\mathbf{D}\delta\phi
- \mathcal{W}[\mathbf{D}\Delta^{(t)}(\mathbf{k})]\|_2^2,
\label{eq:increment_ls_main}
\end{equation}
and the corresponding first-order optimality condition is
\begin{equation}
\mathbf{L}\delta\phi^{(t+1)}(\mathbf{k})
= \mathbf{D}^{\top}\mathcal{W}[\mathbf{D}\Delta^{(t)}(\mathbf{k})].
\label{eq:increment_poisson_main}
\end{equation}

To overcome the limitations of $\ell_2$ formulation (Eq. \eqref{eq:poisson}), we retain the integer offset field $\mathbf{k}$ as an explicit optimization variable rather than eliminating it through the Itoh relation. We assume that $\mathbf{k}$ is piecewise constant within anatomically homogeneous tissue, with transitions occurring at boundaries. The resulting mixed discrete-continuous problem, with ($\delta\phi$ , $\mathbf{k}$) motivates the bi-level decomposition described below.

\subsection{Bi-Level Optimization Framework}
\label{sec:bilevel}

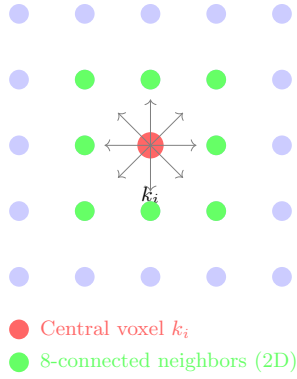
\begin{figure}[t]
\centering
\resizebox{0.45\linewidth}{!}{%
\begin{tikzpicture}
      \begin{scope}
        \foreach \x in {0,1,2,3,4} {
            \foreach \y in {0,1,2,3,4} {
                \fill[blue!20] (\x,\y) circle (0.15);
            }
        }

        \fill[red!60] (2,2) circle (0.2);
        \node[below] at (2,1.5) {\small $k_i$};

        \fill[green!60] (2,3) circle (0.15);
        \fill[green!60] (2,1) circle (0.15);
        \fill[green!60] (1,2) circle (0.15);
        \fill[green!60] (3,2) circle (0.15);
        \fill[green!60] (1,3) circle (0.15);
        \fill[green!60] (3,3) circle (0.15);
        \fill[green!60] (1,1) circle (0.15);
        \fill[green!60] (3,1) circle (0.15);

        \foreach \angle in {0,45,90,135,180,225,270,315} {
            \draw[->,gray] (2,2) -- ($(2,2)+(\angle:0.7)$);
        }
    \end{scope}

    \begin{scope}[shift={(0,-0.8)}]
        \fill[red!60] (0,0) circle (0.15) node[right=2mm] {\small Central voxel $k_i$};
        \fill[green!60] (0,-0.5) circle (0.15) node[right=2mm] {\small 8-connected neighbors (2D)};
    \end{scope}
\end{tikzpicture}%
}
\caption{Connectivity pattern of the regularization term $\|\mathbf{T}\mathbf{k}\|_1$ in~\eqref{eq:upper_objective} for the 2D implementation. Each pixel $k_i$ is coupled to its eight spatial neighbors, with edge weights $w_{ij}$ derived from the quality map and the binary mask.}
\label{fig:regularization_pattern}
\end{figure}

Bi-level optimization formulates hierarchical problems in which one optimization task constrains another~\cite{colson2007overview,dempe2002foundations}. The canonical form involves an upper-level variable $\mathbf{x}_u$ and a lower-level variable $\mathbf{x}_l$:
\begin{subequations}
\begin{align}
\min_{\mathbf{x}_u} \quad & F(\mathbf{x}_u,\mathbf{x}_l) \\
\text{subject to} \quad &
\mathbf{x}_l \in \arg\min_{\mathbf{x}_l} f(\mathbf{x}_u,\mathbf{x}_l).
\end{align}
\label{eq:bilevel_general}
\end{subequations}
For the proposed phase-unwrapping method, the upper-level variable is the integer offset field $\mathbf{k}$ and the lower-level variable is the continuous phase increment $\delta\phi$. At iteration $t$, the bi-level problem is
\begin{subequations}
\begin{align}
&\min_{\mathbf{k} \in \mathcal{K}} \quad 
\|\phi^{(t)} + \delta\phi^{(t+1)}(\mathbf{k})
- \phi_w - 2\pi\mathbf{k}\|_2^2  + \lambda_t\|\mathbf{T}\mathbf{k}\|_1,
\label{eq:upper_objective}
\\
&\text{subject to} \quad 
\delta\phi^{(t+1)}(\mathbf{k})
\in \arg\min_{\delta\phi \in \R^N} \|\mathbf{D}\delta\phi
- \mathcal{W}[\mathbf{D}\Delta^{(t)}(\mathbf{k})]\|_2^2.
\label{eq:lower_objective}
\end{align}
\label{eq:bilevel_ilpu}
\end{subequations}
Here, $\mathcal{K}=\{\mathbf{k} \in \Z^N: k_{\min} \leq k_i \leq k_{\max}\}$ is the admissible set of integer offset maps, with bounds controlling the maximum number of wraps. The first term in~\eqref{eq:upper_objective} enforces consistency with the phase-jump model, whereas $\|\mathbf{T}\mathbf{k}\|_1$ promotes a piecewise-constant offset map. For the 2D implementation,
\begin{equation}
\|\mathbf{T}\mathbf{k}\|_1
= \sum_i \sum_{j \in \mathcal{E}_8(i)} w_{ij}|k_i-k_j|,
\label{eq:graph_tv_main}
\end{equation}
where $\mathcal{E}_8(i)$ denotes the eight-connected neighborhood of pixel $i$, as shown in Fig.~\ref{fig:regularization_pattern}. The edge weights are defined as
\begin{equation}
w_{ij} = \min(w_i,w_j)M_iM_j,
\label{eq:graph_weights_mask}
\end{equation}
where $w_i$ is the normalized quality value and $M_i$ is the binary mask (Appendix \ref{app:mask_generation}) value at pixel $i$. As a result, two neighboring offsets are coupled only when both pixels belong to the valid support of the mask.
The regularizer in Eq.~\eqref{eq:graph_tv_main} encourages a piecewise-constant offset map within valid tissue regions without forcing continuity outside the mask.

The bi-level problem~\eqref{eq:bilevel_ilpu} is solved approximately by alternating between two sub-problems. First, the current offset map $\mathbf{k}^{(t)}$ is fixed and the lower-level phase increment is obtained from one DCT-based Poisson solution. 
Second, the updated continuous phase is fixed and the upper-level offset map is refined by a restricted local search. The lower- and upper-level updates are detailed in Appendices~\ref{app:lower_level} and~\ref{app:upper_level}, respectively.

\subsection{Convergence and Complexity Analysis}

\subsubsection{Convergence Properties}

For a fixed offset map $\mathbf{k}^{(t)}$, the lower-level phase update~\eqref{eq:increment_poisson_main} is obtained from one DCT-based Poisson solve. The upper-level update~\eqref{eq:greedy_selection} approximately minimizes the objective in~\eqref{eq:upper_objective}, with each pixel evaluating three candidates while its neighboring offsets are held at their previous values. The full procedure can be viewed as an alternating approximation of the bi-level model in~\eqref{eq:bilevel_ilpu}. Iteration stops once any of three conditions is met. The first is stabilization of the offset map on the tissue mask, where $\mathbf{k}^{(t+1)} \approx \mathbf{k}^{(t)}$. The second is a masked residual that falls below the prescribed threshold $\varepsilon$. The third is reaching the maximum iteration count $T_{\max}$. In practice, convergence occurs within five iterations for most datasets and within twenty for the more challenging cases.

\subsubsection{Computational Complexity}

DCT and inverse DCT (IDCT) operations dominate at $O(N \log N)$, while the Laplacian source computation, greedy offset search over the voxel-neighborhood graph, and quality-guided smoothing are each $O(N)$. Total complexity per outer iteration is therefore $O(N \log N)$. If $T$ denotes the number of outer iterations, the overall algorithm runs in $O(TN \log N)$ with $T \approx 2$--$20$, reducing to $O(N \log N)$ in practice. In summary, Algorithm~\ref{alg:ilpu} presents the complete ILPU pseudocode.

\begin{algorithm}[t]
\caption{ILPU algorithm}
\label{alg:ilpu}
\begin{algorithmic}[1]
\STATE \textbf{Input:} Wrapped phase $\phi_w$, magnitude/quality map $\mathbf{m}$, parameters $\lambda_0$, $\lambda_{\min}$, $\varepsilon$, $T_{\max}$
\STATE \textbf{Output:} Unwrapped phase $\phi^{(t+1)}$
\STATE \textbf{Initialization:}
\STATE Compute $\mathbf{w}$ and tissue mask $\mathbf{M}$; precompute $\boldsymbol{\Lambda}$
\STATE Compute initial $\phi^{(1)}$ by Laplacian unwrapping of $\phi_w$
\STATE Compute initial $\mathbf{k}^{(1)}$ from the phase-offset estimate
\REPEAT
    \STATE \textbf{Sub-Problem 1, lower level:}
    \STATE \quad Form $\Delta^{(t)}(\mathbf{k}^{(t)})$ via~\eqref{eq:regularized_residual_main}
    \STATE \quad Solve $\delta\phi^{(t+1)}$ via~\eqref{eq:dct_solution} 
    \STATE \quad Update $\phi^{(t+1)} \gets \phi^{(t)} + \delta\phi^{(t+1)}$
    \STATE \textbf{Sub-Problem 2, upper level:}
    \STATE \quad Update $\lambda_t$ via~\eqref{eq:adaptive_lambda}
    \STATE \quad Compute $\mathbf{k}^{(t+1)}$ via~\eqref{eq:greedy_selection}
\UNTIL{$\mathbf{k}^{(t+1)} \approx \mathbf{k}^{(t)}$ on $\mathbf{M}$, or $\|\mathbf{M}\odot(\phi^{(t+1)} - \phi_w - 2\pi \mathbf{k}^{(t+1)})\|_2 < \varepsilon$, or $t \geq T_{\max}$}
\RETURN $\phi^{(t+1)}$
\end{algorithmic}
\end{algorithm}

\section{Experimental Results}
\label{sec:results}

\subsection{Dataset and Experimental Setup}

We evaluated ILPU against established phase-unwrapping methods using brain 2D and 3D MRE data of healthy volunteers. MRE provides a useful validation setting for phase unwrapping because the externally induced tissue motion is harmonic in time with a known drive frequency, which provides an additional temporal consistency criterion for manual unwrapping. Written informed consent was obtained from all participants before data acquisition.

The 2D multislice dataset comprised 3,800 brain phase images from five healthy volunteers, acquired during multifrequency MRE acquisitions at 3 Tesla. Images were obtained using a spin-echo echo planar imaging (EPI) sequence with TE/TR = 30/3600 ms, matrix size $136 \times 136 \times 9$ with a 1 mm slice gap and resolution $1.6 \times 1.6 \times 2$ mm$^3$. Mechanical vibrations were applied at frequencies ranging from 5 to 50 Hz (12 frequencies) and encoded using flow-compensated motion-encoding gradients (MEGs) to consecutively measure tissue displacements along three orthogonal spatial directions: left-right (MEG 1), anterior-posterior (MEG 2), and cranial-caudal (MEG 3). 

The 3D dataset was acquired using a spin-echo EPI-based 3D variant. The acquisition parameters were TE/TR = 79/700 ms, resolution $1 \times 1 \times 1$ mm$^3$, 120 slices phase and partition partial Fourier factors of 6/8 and 5/8, acceleration factors of 2 in the phase direction and 3 in the partition direction, 32 time steps, vibration frequency 16 Hz and matrix size $224 \times 224 \times 120$. For this dataset, we used slices $3$-$114$ from MEG 2.

Each phase image underwent manual unwrapping by an expert observer, verified for both spatial continuity and temporal consistency with the expected harmonic motion pattern. These manually unwrapped images served as the ground truth for quantitative error assessment. Quality maps were derived from magnitude images using normalized intensity values $\mathbf{w} = \mathbf{m}/\max(\mathbf{m})$, where $\mathbf{m}$ represents the magnitude signal. The tissue masks were generated from the same magnitude images using the procedure described in Appendix~\ref{app:mask_generation}.

We compared ILPU against three representative methods spanning different algorithmic paradigms, including (1) standard Laplace unwrapping~\cite{ghiglia1994robust,schofield2003fast}, (2) Flynn~\cite{flynn1997minimum}, and (3) SEGUE~\cite{karsa2019segue}. All methods were applied to the same input data, consisting of wrapped phase images. However, unlike standard Laplace, Flynn and ILPU, SEGUE requires a brain mask as an additional input. 

\subsection{Error Metrics}

We employed two complementary error metrics to assess unwrapping accuracy:

\begin{itemize}
    \item \textbf{Residual error} quantifies the normalized deviation from ground truth across all valid pixels within the tissue mask:
\begin{equation}
E_{\text{residual}} = 100 \times
\frac{\|\mathbf{M}\odot(\phi^{(t+1)} - \phi_{g})\|_2}
{\|\mathbf{M}\odot\phi_{g}\|_2},
\end{equation}
where $\phi_{g}$ is the manually unwrapped reference and $\odot$ denotes element-wise multiplication. This metric directly measures solution accuracy relative to ground truth, reported as a percentage.
 \item \textbf{Rewrap error} 
 computed by re-wrapping the unwrapped result and comparing to the original measurements using a wrapped phase difference: 
 \begin{equation}
E_{\text{rewrap}} = \|\mathbf{M}\odot\mathcal{W}[\phi_{w} - \mathcal{W}(\phi^{(t+1)})]\|_2.
\end{equation}
 This metric detects violations of the fundamental unwrapping constraint that rewrapping the solution must recover the original wrapped phase, serving as an indicator of internal consistency independent of ground truth.
\end{itemize}

For both metrics, smaller values indicate better unwrapping. 
We define perfect reconstruction as cases where the residual error falls below $10^{-6}$\%, a threshold chosen to capture both exact numerical zeros and errors that are negligible relative to a single $2\pi$ pixel deviation from the reference.

\begin{figure}[t]
    \centering
    \includegraphics[width=0.85\linewidth]{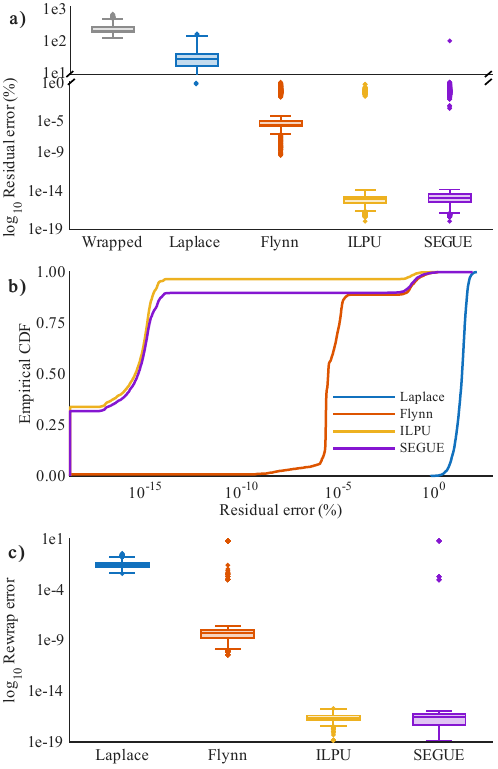}
    \caption{Comprehensive performance comparison across 3,800 brain MRE phase images. (a) Box plot of the residual error (\%) for the wrapped phase and each unwrapping method. A broken axis is used to separate the large errors from the wrapped and Laplace errors from the low‑error methods. (b) ECDF of residual error (\%) on a log-scaled x-axis, emphasizing accumulation near zero. (c) Box plot of rewrap error.}
    \label{fig:box_plots}
\end{figure}

{\begin{figure*}[ht!]
    \centering
    \includegraphics[width=\linewidth]{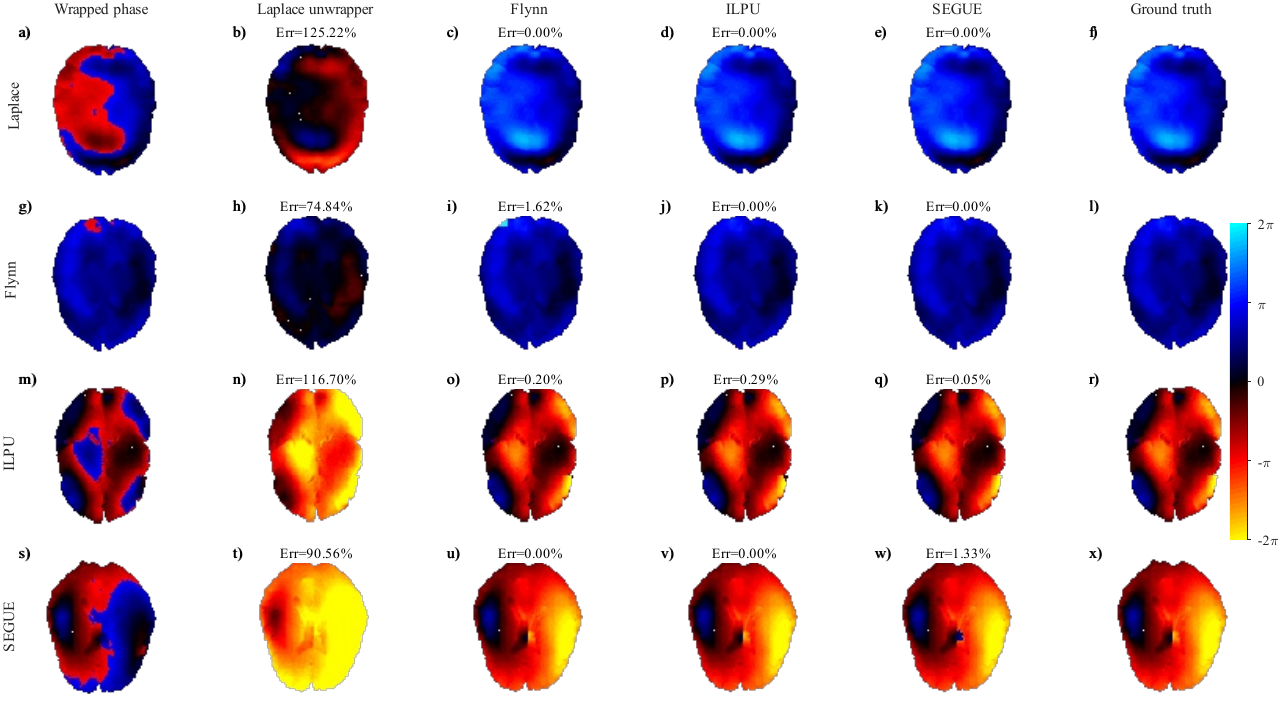}
    \caption{Worst-case visual comparison of phase-unwrapping methods. Each row shows a slice selected as a worst case for one of the compared methods, and the same slice is displayed across all methods together with the wrapped input and ground truth for direct comparison. Within each row, panels share a symmetric color scale derived from the ground truth image.}
    \label{fig:worst_cases}
\end{figure*}}

\subsection{Comprehensive Performance Analysis}
Fig.~\ref{fig:box_plots} provides three complementary views of performance across all 3,800 test images. Fig.~\ref{fig:box_plots}a uses a broken-axis box plot of residual error to separate the much larger Laplace errors from the low-error methods while preserving comparison among Flynn, ILPU, and SEGUE. Fig.~\ref{fig:box_plots}b shows the Empirical Cumulative Distribution Function (ECDF) of residual error, highlighting how rapidly each method accumulates near zero error. Fig.~\ref{fig:box_plots}c reports rewrap error and confirms that the low-error methods remain stable after rewrapping, although residual error provides the clearer separation among the best-performing methods. 

Among the low-error methods, ILPU shows the most favorable residual-error distribution (Fig.~\ref{fig:box_plots}a). Its median residual error is at machine precision ($2.66\times10^{-16}$), and 95\% of cases remain essentially error-free ($4.56\times10^{-15}$). The ECDF in Fig.~\ref{fig:box_plots}b rises most rapidly near zero for ILPU, indicating the largest proportion of near-perfect reconstructions. Compared with Flynn and SEGUE, the main advantage of ILPU appears in the upper tail: the 99th-percentile residual error is 0.10\% for ILPU, versus 0.41\% for Flynn and 0.31\% for SEGUE. Table~\ref{tab:residual_error_summary} summarizes these trends. Only 3.42\% of ILPU cases exceed $10^{-6}\%$ residual error, 1.00\% exceed 0.1\%, and none of them exceed 1\%, all lower than Flynn and SEGUE. 

Over the full 3,800-image dataset, ILPU attains the lowest rewrap error (Fig.~\ref{fig:box_plots}c), with a mean of $2.71\times10^{-18}$ and a maximum of $1.84\times10^{-16}$, compared with a mean of $0.08$ and a maximum of $6.28$ for Flynn, a mean of $0.02$ and a maximum of $6.28$ for SEGUE, and a mean of $0.04$ and a maximum of $0.38$ for Laplace. The larger maxima for Flynn and SEGUE reflect occasional rewrap inconsistencies, whereas ILPU and Laplace remain free of such jumps.

In terms of runtime, the Laplace method required approximately 1.5 ms per image, confirming the computational efficiency of its fast DCT-based solution. ILPU required approximately 23 ms per image, compared with approximately 112 ms per image for Flynn. We do not report a runtime for SEGUE, since its source code was not available to us and the timings obtained from the provided executable were excessively high and appeared unreliable. All runtimes were obtained by processing the full set of 3,800 images on a standard desktop computer (Intel Core i5-12500, 3.0 GHz, 6 cores, 16 GB RAM) and dividing the total time by the number of images.

Together, these results show that ILPU combines near-zero residual error in most images with a smaller tail of failures than the competing low-error methods in a reasonable time.

\begin{table}[!t]
    \centering
    \caption{Summary of residual-error distributions across all 3,800 test images. All values are reported in percent. Smaller values and lower exceedance rates indicate better performance.}
    \label{tab:residual_error_summary}
    \scriptsize
    \setlength{\tabcolsep}{3pt}
    \begin{tabular}{lcccccc}
        \hline
        Method & Median & 95th pct & 99th pct & $>10^{-6}$ & $>0.1$ & $>1$ \\
        \hline
        Laplace & 84.08 & 188.38 & 279.03 & 100 & 100 & 100 \\
        Flynn & $3.00\times10^{-6}$ & 0.09 & 0.41 & 94.24 & 4.79 & 0.21 \\
        \textbf{ILPU} & \textbf{$2.66\times10^{-16}$} & \textbf{$4.56\times10^{-15}$} & \textbf{0.10} & \textbf{3.42} & \textbf{1.00} & \textbf{0.00} \\
        SEGUE & $3.84\times10^{-16}$ & 0.08 & 0.31 & 10.24 & 4.03 & 0.21 \\
        \hline
    \end{tabular}
\end{table}

\subsection{Worst-Case Visual Comparison}
Fig. ~\ref{fig:worst_cases} presents representative worst-case slices selected separately for the compared methods. For each selected slice, the outputs of all methods are shown alongside the wrapped input and the ground truth, allowing direct visual comparison on exactly the same challenging example.

The worst-case comparison is consistent with the quantitative results in Fig.~\ref{fig:box_plots} and Table~\ref{tab:residual_error_summary}. On slices selected by the competing methods, the Laplace solution shows large spatially extended errors, whereas ILPU more closely follows the ground truth. On the ILPU worst-case slice, the remaining discrepancies are visually confined to limited low-SNR or boundary regions rather than spreading across the entire brain. Overall, Fig.~\ref{fig:worst_cases} indicates that when ILPU fails, its errors tend to remain localized and structurally modest compared with the broader deviations observed for Laplace and, in the most difficult cases, Flynn.

\subsection{Spatial Error Distribution}
\begin{figure}[h!t]
    \centering
    \includegraphics[width=0.85\linewidth]{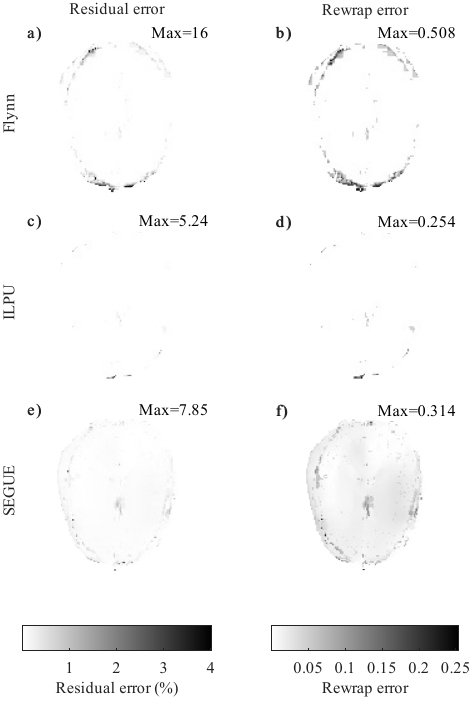}
    \caption{Spatially aggregated error maps across 3,800 images after excluding voxels supported by fewer than 10 images. Residual error (\%) for (a) Flynn, (c) ILPU, and (e) SEGUE, and the corresponding rewrap error for (b) Flynn, (d) ILPU, and (f) SEGUE. Shared grayscale color scales are used within each column.}
    \label{fig:spatial_error}
\end{figure}

Fig. ~\ref{fig:spatial_error} summarizes spatially aggregated residual and rewrap errors for Flynn, ILPU, and SEGUE after excluding voxels supported by fewer than 10 images. Across the retained mask, voxel support ranged from 16 to 3,800 images, with a median of 3,800.

Panels (a), (c), and (e) show the mean residual-error maps. ILPU exhibits the lowest residual error across most of the brain, with a mean of 0.01\%, a median of $1.92 \times 10^{-15}$\%, a 95th percentile of $2.92 \times 10^{-15}$\%, and a maximum of 5.24\%. In comparison, Flynn yields a mean residual error of 0.05\% and a maximum of 16.0\%, while SEGUE yields a mean of 0.07\% and a maximum of 7.85\%.

Panels (b), (d), and (f) show the corresponding rewrap-error maps and preserve the same ordering. ILPU again achieves the lowest values, with mean rewrap error $6.84 \times 10^{-4}$, median $9.38 \times 10^{-17}$, 95th percentile $1.22 \times 10^{-16}$, and maximum 0.25, compared with 0.01 and 0.51 for Flynn and 0.02 and 0.31 for SEGUE. These aggregated maps indicate that ILPU not only reduces overall error but also limits the spatial extent and severity of localized failures.

\subsection{Robustness to Noise and Phase Dynamic Range}
\begin{figure}[h!t]
    \centering
    \includegraphics[width=\linewidth]{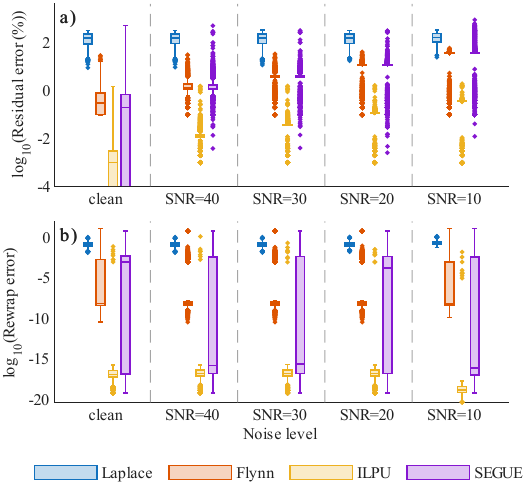}
    \caption{Robustness analysis for the fixed $3 \times$ phase-scaling condition under additive Gaussian noise. Within each noise group (clean, 40 dB, 30 dB, 20 dB, 10 dB SNR), the four methods are Laplace, Flynn, ILPU, and SEGUE. Panel (a) shows $\log_{10}$ residual error (\%), and panel (b) shows $\log_{10}$ rewrap error.}
    \label{fig:robustness}
\end{figure}

To probe robustness under more demanding wrapping conditions, we fixed the phase scaling to 3$\times$ and varied the noise level from clean data to 10 dB SNR. Fig.~\ref{fig:robustness}a summarizes the resulting residual-error distributions for Laplace, Flynn, ILPU, and SEGUE, while Fig.~\ref{fig:robustness}b shows the corresponding rewrap errors.

Across all noise levels, ILPU preserves the lowest residual errors. Its median residual error increased only from $9.74\times10^{-4}$\% in the clean case to 0.36\% at 10 dB SNR, whereas Flynn increased from 0.30\% to 36.0\%, SEGUE from 0.20\% to 36.4\%, and Laplace remained above 157\% throughout. The separation is also evident in the upper tail, where at 10 dB SNR the 95th-percentile residual error was 0.38\% for ILPU, compared with 38.5\% for Flynn, 106.6\% for SEGUE, and 301.5\% for Laplace. Thus, even under severe noise after 3$\times$ phase scaling, ILPU shows graceful degradation rather than catastrophic failure.

Panel~(b) confirms that ILPU also maintains the most controlled rewrap behavior. Its median rewrap error stayed at machine precision across all noise levels, and even at 10 dB SNR its maximum rewrap error was only 0.02. By contrast, Flynn and SEGUE showed much broader rewrap tails, with large 99th-percentile values and maxima of 12.57 and 13.27, respectively, at 10 dB SNR. Laplace retained comparatively modest rewrap error but persistently large residual error, consistent with unresolved integer-offset ambiguities.


\begin{figure*}[h!t]
    \centering
    \includegraphics[width=\linewidth]{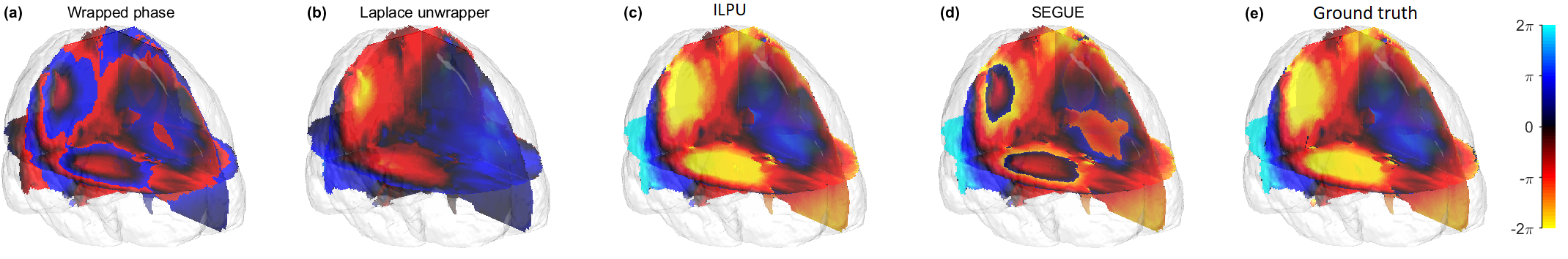}
    \caption{Representative 3D phase-unwrapping example shown with orthogonal cut planes through a cropped brain volume. Panels show (a) wrapped phase, (b) Laplace unwrapper, (c) ILPU, (d) SEGUE, and (e) ground truth. All panels use the same color limits for direct visual comparison. Relative errors with respect to the ground truth were 81.02\% for Laplace, 2.12\% for ILPU, and 67.59\% for SEGUE.}
    \label{fig:3d_phase}
\end{figure*}

These robustness results indicate that the advantage of ILPU persists across the full noise range. Compared with Flynn and SEGUE, the gap widens as SNR decreases, as at 20 dB SNR the median residual error is 0.11\% for ILPU versus 11.5\% for Flynn and 11.4\% for SEGUE, and at 10 dB SNR the corresponding values are 0.36\%, 36.0\%, and 36.4\%. At the same time, the rewrap panel shows that ILPU avoids the large outliers still present for Flynn and SEGUE. These results confirm that the bi-level optimization framework remains stable even when increased wrap density is combined with substantial noise.

\subsection{Three-Dimensional Volume Example}

To illustrate performance on volumetric data, Fig.~\ref{fig:3d_phase} shows a representative cropped 3D phase volume visualized with orthogonal slices through the wrapped input, Laplace solution, ILPU, SEGUE, and the ground truth. All panels use the same color scale, allowing direct comparison of phase structure across methods.

ILPU follows the ground truth volume closely, whereas both Laplace and SEGUE show substantial deviations. The relative error with respect to the ground truth was 2.12\% for ILPU, compared with 67.59\% for SEGUE and 81.02\% for Laplace. The 32-fold reduction in error relative to SEGUE and the 38-fold reduction relative to Laplace demonstrate that the bi-level framework retains its advantage in 3D, where region-growing methods are known to accumulate errors across the larger search space of volumetric data.

\section{Discussion}
\label{sec:discussion}

ILPU achieves superior accuracy through bi-level optimization that addresses limitations of purely transform-based or global optimization approaches. Fig.~\ref{fig:box_plots} and Table~\ref{tab:residual_error_summary} show that this advantage is not only apparent by low mean values but also in the smaller tails of the error distribution. The ECDF shows that ILPU accumulates near-zero errors most rapidly, and the rewrap-error comparison remains consistent with this overall ranking. The worst-case visual comparison in Fig.~\ref{fig:worst_cases} supports the same conclusion qualitatively, where even on particularly difficult slices ILPU remains closer to the ground truth, and its residual errors tend to stay localized near challenging boundary or low-SNR regions rather than propagating broadly through the image. Fig.~\ref{fig:spatial_error} further shows that this advantage persists in the spatially aggregated maps. In comparison with Flynn and SEGUE, ILPU exhibits the lowest mean residual and rewrap errors together with the smallest maxima, indicating that its failures are both less frequent and less spatially pronounced.

The Laplace method provides computational efficiency but fails to recover correct integer offsets, as evidenced by the large residual-error distribution in Fig.~\ref{fig:box_plots}a and its 84.08 median residual error. This arises because the Poisson equation removes DC and low-frequency components, leaving offset ambiguities unresolved. In MRE, this bias is particularly consequential: an oversmoothed offset phase field produces a systematically biased wave curvature, which inversion algorithms interpret as a shorter mechanical wavenumber, causing underestimation of tissue stiffness, particularly at lower frequencies and larger wave numbers.

Flynn employs global optimization to minimize weighted discontinuities. Although widely considered as an accuracy benchmark, Flynn is rarely used in MRE post-processing because its computational cost scales poorly with image count \cite{gdeisat2018performance}. Processing thousands of phase images per subject at the required temporal resolution is impractical without dedicated computing infrastructure. ILPU is designed to bridge the gap between accuracy-first and speed-first methods, exceeding Flynn accuracy while approaching Laplace computational speed through bi-level decomposition that restricts combinatorial search to a three-candidate local neighborhood. Using the $10^{-6}$\% threshold for perfect reconstruction, Flynn achieves only 5.76\% perfect cases on the original dataset compared to 96.58\% for ILPU.  Fig.~\ref{fig:worst_cases} further shows that on the most difficult slices Flynn can produce localized boundary artifacts even when the global structure is largely recovered. The fixed-3$\times$ robustness experiment in Fig.~\ref{fig:robustness} quantitatively demonstrates this limitation. For this test Flynn and SEGUE both deteriorate rapidly as noise increases, reaching median residual errors of 36.0\% and 36.4\% at 10 dB SNR, whereas ILPU remains at 0.36\%. Flynn and SEGUE often preserve low median rewrap error, but their distributions retain pronounced rewrap-error tails and large outliers, pointing to instability where multiple wraps occur between adjacent voxels.

ILPU addresses both limitations through iterative coupling between smooth estimation and discrete refinement. The Laplacian solution provides a high-quality initial estimate. The greedy search then resolves offset ambiguities by evaluating only three candidates per voxel ($\pm 1$ around the Laplacian-based estimate), keeping computational cost low while correcting typical wrap errors. Quality- and mask-weighted regularization weakens coupling in low-quality regions and across mask boundaries while encouraging piecewise-constant offsets within reliable tissue regions.


The computational cost of ILPU, approximately 23 ms per image, represents a favorable trade-off between the high computational speed of Laplace unwrapping, approximately 1.5 ms per image, and the accuracy requirements of clinical applications. For MRE post-processing workflows, where datasets typically contain hundreds to thousands of phase images, the processing time of ILPU remains practical on a standard desktop computer. The roughly 15-fold increase over pure Laplacian unwrapping yields 96.58\% perfect reconstructions against none, a substantial accuracy gain that justifies the additional computation. Relative to Flynn, at approximately 112 ms per image, ILPU achieves a 4.8-fold speedup while maintaining higher accuracy, and unlike Flynn it does not suffer from $2\pi$ constant shifts. ILPU also scales more favorably than global methods for large scale problems, since each iteration remains at $O(N\log N)$ complexity while the search space of global methods grows exponentially with image size.


The ILPU algorithm is formulated so that the Laplacian operator and DCT eigenvalue decomposition generalize directly to three-dimensional volumes at $O(N \log N)$ complexity. The graph-TV regularizer can likewise be extended from the eight-connected 2D graph to a selected voxel-neighborhood graph in 3D. The quantitative evaluation presented in Section~\ref{sec:results} uses 2D images because the majority of MRE reconstruction pipelines currently operate slice-by-slice~\cite{pmid40720631,pmid37705467,pmid35691940}, making 2D unwrapping the most clinically relevant target. Fig.~\ref{fig:3d_phase} demonstrates that the algorithm transfers directly to volumetric data and ILPU achieves a relative error of 2.12\% with respect to the ground truth, compared with 67.59\% for SEGUE and 81.02\% for Laplace. This result is particularly notable because SEGUE, which performs competitively in 2D, degrades substantially in 3D, consistent with the known sensitivity of region-growing methods to error accumulation in volumetric unwrapping. Three-dimensional unwrapping introduces additional challenges including larger memory footprints, longer per-volume processing times, and more complex residue topologies. A systematic volumetric benchmark remains a subject of future work. 


Several limitations warrant discussion. First, ILPU requires quality maps from anatomical images, such as the signal magnitude of MRE images. When magnitude information is unavailable, alternative measures such as phase gradient consistency or local signal stability could provide reliability estimates. Second, the three-candidate neighborhood may prove insufficient for extremely rapid phase variations. However, the fixed-3$\times$ robustness tests demonstrate graceful degradation rather than catastrophic failure even under substantial noise.

Third, the current implementation processes images independently, not exploiting temporal or volumetric correlations. Extension to joint 2D+time and 3D+time optimization will be a useful and straightforward step, since the bi-level framework generalizes naturally to volumetric and temporal domains. Importantly, however, many MRE reconstruction pipelines currently operate slice-by-slice or on 2D phase maps~\cite{pmid40720631,pmid37705467,pmid35691940} and would directly benefit from the proposed algorithm without any further extension.\\

\section{Conclusion}
\label{sec:conclusion}
We presented ILPU, an MRI phase unwrapping algorithm based on bi-level optimization. The method decomposes unwrapping into coupled optimization tasks: a lower level solving an incremental Poisson equation via DCT, and an upper level refining integer offsets through quality-guided spatial regularization using a restricted local search.

The algorithm is formulated for both 2D and 3D data. On 3,800 2D brain MRE images and a $10^{-6}$\% residual-error threshold, ILPU achieved 96.58\% perfect cases versus 89.76\% for SEGUE, 5.76\% for Flynn, and 0\% for Laplace. In 3D, relative errors were 2.12\% for ILPU versus 67.59\% for SEGUE and 81.02\% for Laplace. The algorithm operates at $O(N\log N)$ complexity and is numerically faster than both Flynn and SEGUE, at the cost of additional computation relative to pure Laplacian unwrapping.

Key innovations include: (1) iterative coupling between continuous Laplacian estimation and discrete offset refinement; (2) quality- and mask-weighted eight-connected spatial regularization in the 2D implementation; (3) adaptive regularization scheduling; and (4) restricted three-candidate local search.

\section*{Code and Data Availability}
The ILPU code and the manually unwrapped reference dataset are available from the corresponding author upon reasonable request. Deployment on the BIOQIC web server (\url{https://bioqic-apps.charite.de}) is currently pending.

\appendix
\section{Mask Generation}
\label{app:mask_generation}

The binary mask $\mathbf{M}$ defines the valid tissue support used by the offset regularization and by the stopping criterion. It is generated from the magnitude or quality image before the unwrapping iterations. First, a monogenic surface-symmetry filter is applied to the magnitude image across multiple spatial scales~\cite{felsberg2001monogenic}. This produces a monogenic feature-symmetry map (MFS), whose values are high in coherent tissue regions and near anatomical boundaries, and low in background or unreliable regions.

The initial support is obtained by thresholding the MFS map at a small fraction of its maximum value. Holes inside the detected support are filled, and the mask is then expanded by a few pixels. This final expansion keeps the near-boundary tissue band where phase wrapping is often most difficult, while still preventing unnecessary regularization outside the object.

\section{Lower Level: Laplacian Phase-Increment}
\label{app:lower_level}

For a fixed offset map $\mathbf{k}^{(t)}$, the lower-level step solves the incremental Poisson equation in Eq.~\eqref{eq:increment_poisson_main}, with the regularized residual $\Delta_t(\mathbf{k}^{(t)})$ defined in Eq.~\eqref{eq:regularized_residual_main}. The implementation evaluates the right-hand side through the complex-domain Laplacian source
\begin{equation}
\rho
= \mathrm{Im}\left[
e^{-\mathrm{i}\Delta_t(\mathbf{k}^{(t)})}
\odot
\mathbf{L}e^{\mathrm{i}\Delta_t(\mathbf{k}^{(t)})}
\right].
\label{eq:complex_source}
\end{equation}
Here, $\rho$ is the numerical source term associated with the current regularized residual. The equation solved at the lower level is therefore
\begin{equation}
\mathbf{L}\delta\phi^{(t+1)} = \rho .
\label{eq:poisson_increment_app}
\end{equation}
The increment is computed with the DCT-II Poisson solver~\cite{ghiglia1994robust,schofield2003fast}:
\begin{equation}
\delta\phi^{(t+1)}
= \mathrm{IDCT}\left(
\frac{\mathrm{DCT}(\rho)}{\boldsymbol{\Lambda}}
\right).
\label{eq:dct_solution}
\end{equation}
For a 3D volume, the DCT eigenvalues are
\begin{equation}
\begin{aligned}
\lambda_{u,v,w}
={}& 2\cos\left(\frac{\pi u}{W+1}\right)
+ 2\cos\left(\frac{\pi v}{H+1}\right) \\
& + 2\cos\left(\frac{\pi w}{S+1}\right) - 6,
\end{aligned}
\label{eq:eigenvalues}
\end{equation}
for $u = 1,\ldots,W$, $v = 1,\ldots,H$, and $w = 1,\ldots,S$. 

\section{Upper Level: Integer Offset Refinement}
\label{app:upper_level}

The upper-level objective is defined in Eq.~\eqref{eq:upper_objective}, and the graph-TV term is given in Eqs.~\eqref{eq:graph_tv_main} and~\eqref{eq:graph_weights_mask}. Direct minimization over all integer maps is computationally expensive. The implementation therefore uses a restricted parallel local search. First, define
\begin{equation}
k_{0,i}
= 
\frac{\phi_i^{(t+1)} - \phi_{w,i}}{2\pi}
.
\label{eq:candidates}
\end{equation}
At each pixel $i$, only the three candidates $k_{0,i}+m$, with $m \in \{-1,0,1\}$, are evaluated. The local cost is
\begin{multline}
C_{i,m}
= \bigl(\phi_i^{(t+1)} - \phi_{w,i} - 2\pi(k_{0,i}+m)\bigr)^2 \\
+ \lambda_t \sum_{j \in \mathcal{E}_8(i)}
\min(w_i,w_j)M_iM_j|k_{0,i}+m-k_j^{(t)}|,
\label{eq:local_cost}
\end{multline}
where $\mathcal{E}_8(i)$ denotes the eight-connected neighborhood of pixel $i$ (Fig.~\ref{fig:regularization_pattern}). The neighboring offsets are held fixed at their values from the previous outer iteration. The offset update is
\begin{equation}
k_i^{(t+1)}
= \Pi_{\mathcal{K}}[k_{0,i} + \arg\min_{m \in \{-1,0,1\}} C_{i,m}].
\label{eq:greedy_selection}
\end{equation}
The projection $\Pi_{\mathcal{K}}$ enforces the bounds and also rounding to the nearest integer.

The regularization parameter is updated according to
\begin{equation}
\lambda_{t+1}
= \max\left(\lambda_{\min},\lambda_t\frac{r_t}{r_0}\right),
\qquad
r_t
= \frac{\|\Delta_t(\mathbf{k}^{(t)})\|_2^2}{\|\phi_w + 2\pi\mathbf{k}^{(t)}\|_2^2},
\label{eq:adaptive_lambda}
\end{equation}
where $r_0$ is the residual value at the first refinement step and $\lambda_{\min}$ prevents excessive reduction of the regularization weight.

\bibliographystyle{unsrt}
\bibliography{references}

\end{document}